\documentclass{amsart}

\usepackage[T1]{fontenc}
\usepackage{ae,aecompl}

\usepackage{amssymb,amscd}
\usepackage{url}
\usepackage{amsthm,amsfonts,amsmath,amsxtra}
\usepackage{graphicx}
\usepackage{tikz} 

\numberwithin{equation}{section}
\numberwithin{table}{section}
\numberwithin{figure}{section}

\newtheoremstyle{bold}
{.5\baselineskip}{.5\baselineskip}{\itshape}{}{\bfseries}{.}{.5em}{}
\newtheoremstyle{shy}
{.5\baselineskip}{.5\baselineskip}{}{}{\bfseries}{.}{.5em}{}

\makeatletter
\def\@captionfont{\small}

\def\section{\@startsection{section}{1}%
  \z@{.9\linespacing\@plus\linespacing}{.5\linespacing}%
  {\large\bfseries\boldmath\centering}}

\def\subsection{\@startsection{subsection}{2}%
  \z@{.7\linespacing\@plus\linespacing}{.5\linespacing}%
  {\normalfont\scshape\centering}}

\def\theindex{\@restonecoltrue\if@twocolumn\@restonecolfalse\fi
  \columnseprule\z@ \columnsep 35\p@
  \@indextitlestyle
  \thispagestyle{empty}%
  \let\item\@idxitem
  \parindent\z@  \parskip\z@\@plus.3\p@\relax
  \raggedright
  \hyphenpenalty\@M
  \footnotesize}

\renewcommand{\@bibtitlestyle}{%
  \@xp\section\@xp*\@xp{\bibname}%
}
\renewcommand{\tocchapter}[3]{%
  \indentlabel{\@ifnotempty{#2}{\ignorespaces#1 #2.\quad}}#3}

\renewcommand{\tocsection}[3]{%
  \indentlabel{\@ifnotempty{#2}{\makebox[3.2em][l]{\ignorespaces#1 #2.}}}#3}

\makeatother

\renewcommand{\bibname}{References}
\renewcommand{\ge}{\geqslant}
\renewcommand{\geq}{\geqslant}
\renewcommand{\le}{\leqslant}
\renewcommand{\leq}{\leqslant}

\theoremstyle{bold}
\newtheorem{theorem}{Theorem}[section]

\theoremstyle{shy}

\newtheorem{remark}[theorem]{Remark}

\newcommand{\cA}{\mathcal{A}}

\newcommand{\cC}{\mathcal{C}}

\newcommand{\EE}{\mathbb{E}}

\newcommand{\NN}{\mathbb{N}}

\newcommand{\PP}{\mathbb{P}\ts}

\newcommand{\RR}{\mathbb{R}}

\newcommand{\ZZ}{\mathbb{Z}\ts}

\newcommand{\ts}{\hspace{0.5pt}}

\title[Ancestral lineages in spatial population models]{Ancestral lineages in spatial population models with local regulation}
\author{Matthias~Birkner and Nina~Gantert}
\date{\today}
\makeindex

\begin{document}

\begin{abstract}
We give a short overview on our work on ancestral lineages in spatial population models with local regulation. 
We explain how an ancestral lineage can be interpreted as a 
random walk in a dynamic random environment. 
Defining regeneration times allows to prove central limit theorems for such walks. We also consider several ancestral lineages in the same population and show for one prototypical example that in one dimension 
the corresponding system of coalescing walks converges to the Brownian web.
\end{abstract}

\maketitle

\newcommand{\citet}[1]{[{\tt #1}]}

\section{Introduction}

Many natural populations live in a spatially extended -- and often
essentially two-dimensional -- habitat, with a range that is much
larger than the typical distance that any individual may travel during
its lifetime. When different genetic types are considered, this can
lead to a local differentiation of types that violates the assumptions
of panmixia. Furthermore, as a result of the interaction of
individuals with their environment -- which may be influenced by the
population itself and additionally by other, competing species or by
external events -- local population sizes often fluctuate in time, and
these fluctuations may be described using random fields.
Understanding the evolution of populations with spatial structure is
an interesting problem, and mathematical, individual-based models can
help to understand how spatial structure modifies the action of other
evolutionary forces such as genetic drift or selection.  \smallskip

It is natural to translate the question of the spatial distribution of types into one
about the spatial embedding of \emph{genealogies} by analysing the
space-time history of sampled individuals and their ancestral
lines. In order to make the latter mathematically tractable, 
a customary approach, especially in mathematical population genetics, 
is to impose a discrete grid of `demes' and assume that local
population sizes are constant in time, as in Kimura's stepping stone
model and its relatives \cite{MBNG-WeissKimura:1965,MBNG-Sawyer:1976}.
\index{stepping stone model}%
Then, ancestral lines of sampled individuals are coalescing random walks
\index{random walks!coalescing}%
(with a delay depending on the local population size), and detailed formulas
for quantities of interest like the decay of the probability of
identity by descent or the correlation of type frequencies with
spatial separation are available
\cite{MBNG-Wilkinson-Herbots:2003,MBNG-ZaehleCoxDurrett:2005}.

Arguably, the built-in assumption of fixed local population sizes in
stepping stone models, though allowing the use of powerful
mathematical tools in the ana\-lysis, appears somewhat artificial from
the modelling perspective.  We remark here that the most `obvious' 
attempt at removing this assumption would be to consider populations 
which move and reproduce freely in space without interaction among them, 
i.e.\ systems of critical branching random walks. The assumption of 
criticality, i.e.\ on average one offspring per individual, is a necessary, 
though not sufficient condition for such systems to possess non-trivial 
equilibria. Unfortunately, this attempt is bound to fail, at least in 
spatial dimensions $d=1$ and $d=2$, although the latter is possibly 
the most interesting case from a biological point of view: 
It is well known that in dimensions $1$ and $2$, critical branching random walks `generically' 
exhibit local extinction and if one conditions on non-extinction, the configuration 
forms arbitrarily dense clumps (\cite{MBNG-Kallenberg}, see, e.g., 
Ch.~6.4 in \cite{MBNG-Etheridge:2011} for a discussion). This effect can also not be 
eliminated by density-dependent down-regulation of the branching rate, see \cite{MBNG-BS19}.
\index{branching random walk}%
\smallskip 

Another line of thought, more in the vein of mathematical
{ecology}, aims at remedying the artificial and in principle
undesirable assumption of fixed local population sizes and some
formulations also remove the discretisation of space in the models
discussed above. Here, one models explicitly the stochastic evolution
of the local population size forward in time in a way that takes
`feedback' into account, typically in the sense that an individual in
a crowded region tends to leave on average less offspring than an
individual that happens to be in a sparsely populated region. Such
models were introduced in the biology literature (and analysed with
non-rigorous methods) in
\cite{MBNG-BolkerPacala:97,MBNG-BolkerPacala:99,MBNG-LawDieckmann:02,MBNG-RaghibHillDieckmann:2011}. Several
investigations in the mathematical literature were inspired by these
models and some modifications thereof, see for instance
\cite{MBNG-NeuhauserPacala:1999, MBNG-BartonDepaulisEtheridge:2002,MBNG-Etheridge:2004,MBNG-BlathEtheridgeMeredith:2007,
  MBNG-FournierMeleard:2004,MBNG-HutzenthalerWakolbinger:07,
  MBNG-BD07} for models and results in this direction
(some with discrete, some with continuous space and `masses').
Models from this class can possess non-trivial equilibria in any spatial dimension 
and they can be `enriched' to also include
ancestral information (this is straightforward for discrete-mass
models as in \cite{MBNG-NeuhauserPacala:1999, MBNG-FournierMeleard:2004,MBNG-BD07}, for
continuous mass models, one could approximate with particle systems or 
use `lookdown' constructions as, e.g., in \cite{MBNG-LePardouxWakolbinger:2013,MBNG-VeberWakolbinger:2014}).
Thus the problem of describing the space-time embedding of ancestral lineages of one or several
individuals sampled from certain locations in an equilibrium
population is mathematically well defined. 
It turns out that a single ancestral lineage, corresponding to a sample of size one, 
then forms a random walk in a dynamic random environment which is generated by the 
backward in time history of the entire population. Similarly, the ancestral information for a larger sample 
corresponds to a system of several random walks in the same environment which 
can additionally coalesce when they are in the same location.
In this article, we discuss the behaviour of ancestral lineages in two prototypical examples, 
namely the discrete time contact process in Section~\ref{MBNG-sect:cp} and the logistic branching random walk 
in Section~\ref{MBNG-sect:lbrw}. A key idea in both Sections~\ref{MBNG-sect:cp} and \ref{MBNG-sect:lbrw} will 
be to construct regenerations. 
It turns out that in both cases, ancestral lineages behave similarly 
as random walks on large space-time scales in the sense that they satisfy the law of large numbers 
and a central limit theorem. Thus, broadly speaking, the effect of the fluctuating local population 
sizes manifests itself on large scales only in the variance parameter of the `random walks'. 
This validates the pragmatical approach mentioned above, where one simply replaces the 
true demographic history of the population by one with locally fixed `effective sizes' (and the migration 
by an `effective migration'). 
In Section~\ref{MBNG-sect:notes}, we discuss 
the relation to 
other projects within SPP 1590 and to the (considerable) literature of random walks in random environments.

\section[RW on the oriented percolation cluster]{The contact process and random walk 
on the backbone of the oriented percolation cluster}
\label{MBNG-sect:cp}

We start with a more detailed description of the model forwards in time and then discuss its ancestral lineages.

\subsection{The discrete time contact process}
\label{MBNG-subsect:dtcp}
Let $\omega := \{\omega(x,n):
(x,n)\in \ZZ^d \times \ZZ\}$ be a family of independent Bernoulli random variables
(representing the carrying capacities) with parameter $p \in (0,1]$. 
We call a site $(x,n)$ \emph{inhabitable} (or \emph{open}) if
$\omega(x,n)=1$ and \emph{uninhabitable} (or \emph{closed}) if $\omega(x,n)=0$. We say that there is an
\emph{open path} from $(y,m)$ to $(x,n)$ for $m \le n$ if there is a sequence
$x_m,\dots, x_n$ such that $x_m=y$, $x_n=x$, $\lVert x_k-x_{k-1} \rVert \le 1$ for
$k=m+1, \dots, n$ and $\omega(x_k,k)=1$ for all $k=m,\dots,n$. In this case we
write $(x,m) \to (y,n)$. Here $\lVert \cdot \rVert$ denotes the $\sup$-norm.
The terms open/closed are standard in percolation theory, we use here
inhabitable/uninhabitable to emphasise the population interpretation.
\index{open}\index{closed}\index{inhabitable}\index{uninhabitable}

Given a set $A \subseteq \ZZ^d$ we define the \emph{discrete time contact process}
$(\eta_n^A)_{n \ge m}$ starting at time $m \in \ZZ$ from the set $A$ as
\begin{align*}
  \eta_m^A (y) & ={\bf 1}_{A}(y), \; y \in \ZZ^d, \\
  \intertext{and for $n \ge m$}
  \eta_{n+1}^A(x) & =
  \begin{cases}
    1 & \text{if $\omega(x,n+1)=1$ and $\eta_n^A(y)=1$ for some $y \in \ZZ^d$
      s.t. $\lVert x-y \rVert \le 1$}, \\
    0 & \text{otherwise}.
  \end{cases}
\end{align*}
\index{contact process!in discrete time}%
In other words, $\eta_n^A(y)=1$ if and only if there is an open path from
$(x,m)$ to $(y,n)$ for some $x\in A$ (where we use in this definition the
convention that $\omega(x,m)={\bf 1}_{A}(x)$ while for $k>m$ the $\omega(x,k)$
are i.i.d.\ Bernoulli as above). 
Taking $m=0$, we set
\begin{align} \label{MBNG-def:tauA}
  \tau^{A} :=  
\inf\{n \ge 0: \eta_n^A \equiv 0 \}.
\end{align}

We interpret the process $\eta$ as a population process, where $\eta_n(x)=1$
means that the position $x$ is occupied by an individual in generation $n$. 
Space-time sites 
can be inhabitable (if $\omega(x,n)=1$) or uninhabitable 
(if $\omega(x,n)=0$). The population dynamics is then the following: 
For each $x \in \ZZ^d$ independently, if $\omega(x,n)=1$ and there was at least one individual in the 
neighbourhood of $x$ in the previous generation, i.e.\ 
\begin{align}
\label{MBNG-eq:defancestors}
\parbox{0.85\textwidth}{if $\cA_{x,n} := \{ y \in \ZZ^d : 
\lVert x-y \rVert \le 1\text { and } \eta_{n-1}(y)=1 \} \neq \emptyset$, \\[0.5ex]
then $y$ is picked uniformly from $\cA_{x,n}$ and an offspring of the individual 
at $y$ in generation $n-1$ is placed at space-time site $(x,n)$.
}  
\end{align}
In this case $\eta_n(x)=1$ and \eqref{MBNG-eq:defancestors} defines the ancestral structure 
of the population.
\index{ancestral structure}%
In the other cases, namely if $\omega(x,n)=0$ (site uninhabitable) or if $\cA_{x,n} = \emptyset$ 
(no inhabited neighbours in the previous generation), we have 
$\eta_n(x)=0$, i.e., the site stays vacant. With this interpretation, \eqref{MBNG-def:tauA} 
is the extinction time of a population which starts with all $x \in A$ inhabited. 
\smallskip

Note that the dynamics \eqref{MBNG-eq:defancestors} implicitly
contain a local population regulation: Neighbours compete for
inhabitable sites, so individuals in sparsely populated regions have
on average higher reproductive success. We can visualise this by 
considering a neutral multi-type version, where offspring simply inherit their parent's type 
(discussed in more detail in Remark~\ref{MBNG-rem:multitpye} below). See the example in Figure~\ref{MBNG-fig1}.
\begin{figure}[h]
  \begin{center}
    \begin{tabular}{ccc}
      \raisebox{0.8cm}{\parbox{0.9cm}{$n$\\[0.7cm]$n-1$}} \hspace{-0.8cm}
& \includegraphics[scale=0.83]{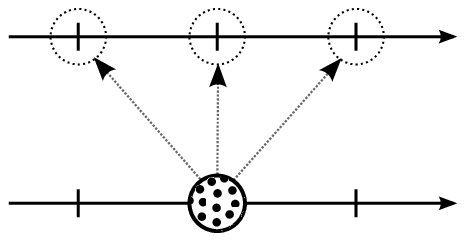} & 
      \includegraphics[scale=0.83]{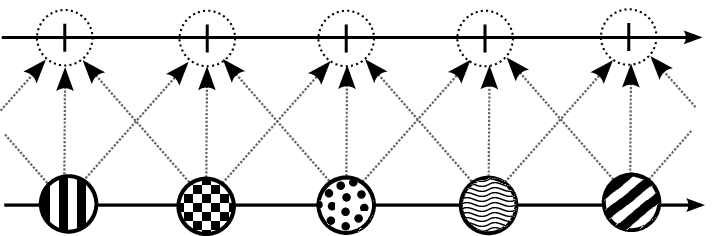} \\[0.0cm]
& \parbox{5cm}{
        \begin{center} expected no.\ of dotted offspring:\\ $3p >
1$
        \end{center} } & 
      \parbox{5cm}{
        \begin{center} expected no.\ of dotted offspring:\\ $3
\frac{1}{3}p = p < 1$
        \end{center} }
    \end{tabular}
  \end{center}
  \caption{Interpretation of $(\eta_n)$ as a locally regulated population model (note that $p > p_c \ge 1/3$ in this case)} 
  \label{MBNG-fig1}
\end{figure}
\index{local population regulation}%

Note that the ancestor is random since we are not given the whole evolution of the system but only its state at time $n$. Compare with the more familiar case of a continuous-time contact process and its graphical representation.  Define the ancestor at time $0$ of an infected site at time $t$ to be the site where the infection came from, following back the graphical representation.
Then, the ancestor at time $0$ of an infected site at time $t$ is determined if we know the graphical representation up to time $t$, but it is random if we only see the configuration of all infected sites at time $t$.

It is well known, see e.g.\ Theorem~1 in \cite{MBNG-GrimmetHiemer:02},
\index{critical value!for oriented percolation}%
that there is a critical value $p_c \in (0,1)$ such that
$\PP(\tau^{\{\boldsymbol 0\}} =\infty )=0$ for $p\le p_c$ and
$\PP(\tau^{\{\mathbf 0\}} =\infty )>0$ for $p> p_c$. Here and in the following, we write
$\mathbf{0} = (0,0,\dots,0) \in \ZZ^d$ for the origin in $d$-dimensional space.

We will only consider the supercritical case $p > p_c$. In this case the
law of $\eta_n^{\ZZ^d}$ converges weakly to the so-called upper invariant measure
which is the unique non-trivial extremal invariant measure of the
discrete-time contact process. By taking $m \to -\infty$ while keeping
$A=\mathbb Z^d$ one obtains the stationary process
\begin{align}
  \label{MBNG-eq:5}
  \eta :=(\eta_n)_{ n \in \ZZ} := (\eta_n^{\ZZ^d})_{n \in \ZZ}.
\end{align}


\subsection{Ancestral lineages}
\label{MBNG-sect:DRWOPC}

We are interested in the behaviour of the `ancestral lineages' of individuals 
in the stationary process $\eta$ from \eqref{MBNG-eq:5}, where 
the behaviour of such a lineage is described by iterating \eqref{MBNG-eq:defancestors}. 
Due to time stationarity, we can focus on ancestral lines of individuals
living at time $0$. It will be notationally convenient to time-reverse the
stationary process $\eta$ and consider the process $\xi := (\xi_n)_{n \in
  \ZZ}$ defined by $\xi_n(x)=1$ if $(x,n) \to \infty$ (i.e.~there is an infinite
directed open path starting at $(x,n)$) and $\xi_n(x)=0$ otherwise. Note that
indeed $\mathcal{L}((\xi_n)_{n\in\ZZ})=\mathcal{L}((\eta_{-n})_{n\in\ZZ})$.
More precisely, due to \eqref{MBNG-eq:5}, $\eta_{-n}(x)= 1$ if and only if there is an infinite directed open backwards
path starting at $(x,-n)$, i.e.~a connection from $-\infty$ to  $(x,-n)$. This is the case if and only if in the time-reversed picture, there is a connection from $(x,n)$ to $\infty$, i.e.~there is an infinite
directed open path starting at $(x,n)$, and this is the case if and only if and only if $\xi_n(x)=1$.
Hence there is a one-to one correspondence of $(\xi_n)_{n\in\ZZ}$ and $(\eta_{-n})_{n\in\ZZ}$ and in particular the two processes have the same law.

We will from now on in this section consider the forwards evolution of $\xi$ as the `positive' 
time direction. \smallskip

On the event $B_0 := \{\xi_0(\mathbf{0})=1\}$ there is an infinite path
starting at $(\mathbf 0,0)$. We define the oriented cluster by
\begin{align*}
  \cC := \{(x,n) \in \ZZ^d \times \ZZ: \xi_n(x)=1\}
\end{align*}
\index{oriented percolation}\index{oriented percolation!cluster}%
(in percolation jargon, this is strictly speaking the `backbone' of
the oriented cluster) and let
\begin{equation}
  \label{MBNG-Udef}
  U(x,n) := \{(y,n+1): \lVert x-y \rVert \le 1\}
\end{equation}
be the neighbourhood of the site $(x,n)$ in the next generation.
One can allow more general finite neighbourhoods in \eqref{MBNG-Udef} with mostly
only notational changes in the proofs, see \cite[Remark~1.4]{MBNG-BCDG13}.
Note however that if $U(x,n)$ is not
symmetric around $x$, the walk will generically have a non-trivial speed.

\begin{figure}[h]
  \includegraphics[width=6.0cm]{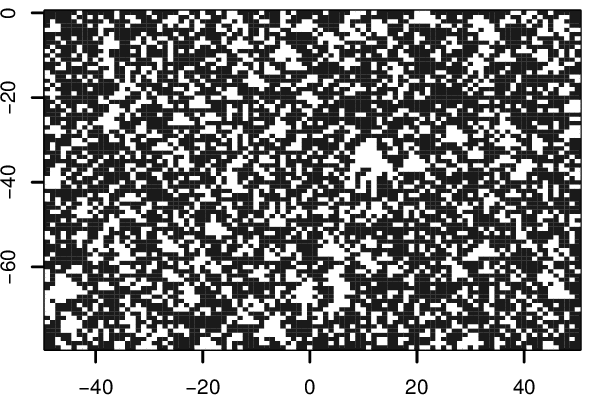}\hspace{0.3cm}
  \includegraphics[width=6.0cm]{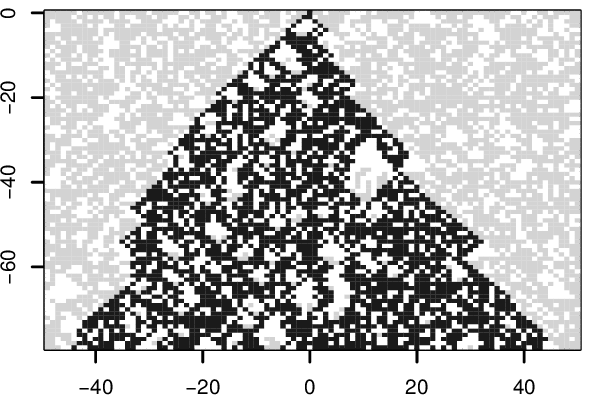}
  \caption{Left: A simulation of the space-time configuration of the stationary
    contact process $\eta$ from \eqref{MBNG-eq:5} with $p=0.68$. Dark 
    sites have $\eta_n(x)=1$. 
    Right: The same configuration with only those sites highlighted in dark 
    which are potential ancestors of the individual at the origin $(\mathbf{0},0)$,
    i.e.\ those sites which the walk $X$ with dynamics \eqref{MBNG-eq:defXdynamics} can reach.
  }
  \label{MBNG-fig:cluster}
\end{figure}

On the event $B_0$ we may define a $\mathbb Z^d$-valued random walk
$X := (X_n)_{n \ge 0}$ starting from $X_0 =\mathbf 0$ with transition
probabilities
\begin{align}
  \label{MBNG-eq:defXdynamics}
  \PP( X_{n+1}=y \mid X_{n}=x, \xi) & =
  \begin{cases}
    {|U(x,n) \cap \cC|}^{-1}  & \text{when }(y,n+1) \in  U(x,n) \cap \cC,\\
    0& \text{otherwise.}
  \end{cases}
\end{align}
\index{ancestral random walk}\index{random walk!on oriented percolation cluster}%
This corresponds to `going backwards' in \eqref{MBNG-eq:defancestors} and we 
interpret $X_n$ as the spatial position of the ancestor $n$ generations ago 
of the individual at the origin today, see also Figure~\ref{MBNG-fig:cluster}.

Note that $(X_n,n)_{n\ge 0}$ is a directed random walk on the
percolation cluster $\cC$, and $X$ can be also viewed as a random walk
in a (dynamical) random environment, where the environment is given by
the process $\xi$. We write $P_\omega$ and $E_\omega$ to denote probabilities 
and expectations when the environment (which is a function of the $\omega$'s) 
is fixed, and write $\PP$ and $\EE$ for the situation when we average 
with respect to both the walk and the environment. In the jargon of random walks 
in random environments, this refers to the `quenched' and the `averaged' or `annealed' case, respectively. 

The main result from \cite{MBNG-BCDG13} is the following theorem on the position $X_n$ of the random walk on the backbone of the oriented percolation cluster at time $n$. The theorem can be interpreted by saying that $X_n$ behaves similarly as a simple random walk: it satisfies a law of large numbers and a central limit theorem. (The case of simple random walk corresponds in our notation to $p=1$.) In other words, the percolation cluster behaves, on large scales, similarly as the full lattice: the effect of the \lq holes\rq\ in the cluster  -- which are clearly visible in the simulation
in Figure~\ref{MBNG-fig:cluster} --
vanishes on large scales.

\begin{theorem}[Law of large numbers, averaged and quenched central limit theorem, 
  {\cite[Theorems~1.1. and 1.3]{MBNG-BCDG13}}]
  \label{MBNG-thmrwopc}
  For any $d \ge 1$ we have
  \begin{align}
    \label{MBNG-eq:LLN}
    P_\omega\Big( \frac1n X_n \to \mathbf{0} \,\Big) = 1
    \quad \text{for $\PP(\,\cdot\, | B_0)$-a.a.\ $\omega$},
  \end{align}
  and for any $f \in C_b(\RR^d)$
  \begin{align}
    \label{MBNG-eq:annealedCLT}
    \EE \Big[ f\left( X_n/\sqrt{n}\, \right) \, \Big| \, B_0 \Big]
    \xrightarrow{n \to
      \infty}   \Phi(f),
  \end{align}
    \begin{equation}
    \label{MBNG-eq:quenchedCLT}
    E_\omega\Big[f\left(X_n/\sqrt{n} \, \right) \,\Big] \xrightarrow{n \to
      \infty} \Phi(f)  \quad    \text{ for  $\PP(\,\cdot\, |
        B_0)$-a.a.~$\omega$,}
  \end{equation}
  where $\Phi(f) := \int f(x)\, \Phi(dx)$ with $\Phi$ a non-trivial
  centred isotropic $d$-dimensional normal law.
  Functional versions of \eqref{MBNG-eq:annealedCLT} and \eqref{MBNG-eq:quenchedCLT} hold as well.
\end{theorem}
A proof sketch is given in Section~\ref{MBNG-subsect:proofBCDG} below.
\index{law of large numbers!for random walk on percolation cluster}%
\index{central limit theorem!for random walk on percolation cluster}%

\begin{remark}
\label{rem:sigma.p}
The covariance matrix of $\Phi$ in \eqref{MBNG-eq:annealedCLT} is
  $\sigma^2$ times the $d$-dimensional identity matrix. It follows from the
  regeneration construction (see Subsection~\ref{MBNG-subsect:proofBCDG} below) that
\begin{equation}
  \sigma^2 = \sigma^2(p) = \frac{\EE\big[ Y_{1,1}^2 \big]}{\EE[\tau_1]} \in (0, \infty)
\end{equation}
where $\tau_1$ is the first regeneration time (see \eqref{MBNG-deftauiYi} below)
of the random walk $X$ and $Y_{1,1}$ is the first coordinate of
$X_{\tau_1}$, the position of the random walk at this regeneration
time.
The behaviour of $\sigma^2(p)$ as $p
\downarrow p_c$ is an interesting open problem that merits
further research.
\end{remark}


\begin{remark}[Consequences for the long-time behaviour of the multi-type 
process] 
\label{MBNG-rem:multitpye}
\index{contact process!multi-type}%
Let us enrich the contact process $(\eta_n)_n$ from Section~\ref{MBNG-subsect:dtcp} 
by including (so-called neutral) types: Say, at time $n=0$, every $\eta_0(x)$ is independently 
assigned a uniformly chosen value from $(0,1)$ and we augment the rule 
\eqref{MBNG-eq:defancestors} by setting $\eta_n(x)=\eta_{n-1}(y) > 0$ if $y \in \cA_{x,n}$ 
was chosen as the ancestor of the individual at site $(x,n)$. Thus, children inherit 
their parent's type (which is $>0$) and we still interpret $\eta_n(x)=0$ as a vacant site. 
As $n\to\infty$, $\eta_n$ will converge in distribution to an equilibrium $\widetilde\eta$ 
of the multi-type dynamics.

It follows from Theorem~\ref{MBNG-thmrwopc} and its proof in \cite{MBNG-BCDG13} 
that any two ancestral lineages will eventually meet in $d \le 2$, but not in $d \geq 3$. 
By `looking backwards in time', this has consequences for $\widetilde\eta$: For any $x,y \in \ZZ$, 
\begin{equation}
\label{MBNG-multitypd12}
\PP\left(\widetilde{\eta}(x) = \widetilde{\eta}(y) \mid \widetilde{\eta}(x)>0, \widetilde{\eta}(y)>0 \right)=1
\end{equation}
in $d=1,2$ and this probability is $<1$ in $d \geq 3$. In fact, for $d \ge 3$ there is $C_d \in (0,\infty)$
such that 
\[
\PP\left(\widetilde{\eta}(x) = \widetilde{\eta}(y) \mid \widetilde{\eta}(x)>0, \widetilde{\eta}(y)>0 \right) 
\sim \frac{C_d}{\lVert x-y\rVert^{d-2}_2} \quad \text{as } \lVert x-y\rVert \to \infty.
\]
These properties are analogous to those of the multi-type stepping stone model. 
\end{remark}

\subsection{Proof ideas: Local construction and regeneration} 
\label{MBNG-subsect:proofBCDG}

A main difficulty in the proof of Theorem~\ref{MBNG-thmrwopc} lies in
the fact that in order to determine $\xi(x,n)$, one has to know the
`whole future' of the environment $\omega $. To overcome this, we
build a trajectory of $X$ using rules that are `local', i.e.~which use
only local $\omega $'s (and some additional local randomness), but not
the $\xi$'s.  We then read off regeneration times from this
construction: These are exactly the times when the locally constructed
trajectory coincides with the true trajectory of $X$, see
\eqref{MBNG-eq:reg.1} below.  This approach is inspired by
\cite{MBNG-Kuczek:1989} and \cite{MBNG-Neuhauser:1992}.

The construction employs some additional randomness: For every
$(x,n) \in \ZZ^d \times \ZZ$ let $\widetilde\omega{(x,n)}$
be a uniformly chosen permutation of $U(x,n)$ ($U(x,n)$ may be written as a vector by ordering the elements according to the lexicographical ordering of the space ccordinate $x$),
independently distributed for all
$(x,n)$'s and independent from the $\omega$'s. 
We denote the whole family of these permutations by $\widetilde\omega$.

For every $(x,n)\in \mathbb Z^d\times \mathbb Z$ let
$\ell(x,n)=\ell_\infty(x,n)$ be the length of the longest directed
open path starting at $(x,n)$; we set $\ell(x,n)=-1$ when $(x,n)$ is
closed. (Recall that a path $(x_0,n), (x_1,n+1), \dots, (x_k,n+k)$ of
length $k$ with $\lVert {x_i-x_{i-1}} \rVert \leq 1$ is open if
$\omega(x_0,n)=\omega(x_1,n+1)=\cdots= \omega(x_k,n+k)=1$. $\ell(x,n)=\infty$ 
means $(x,n) \in \cC$.) For every
$k\in \NN_0$ 
let $\ell_k(x,n) := \ell(x,n)\wedge k$ be
the length of the longest directed open path of length at most $k$
starting from $(x,n)$. Observe that $\ell_k(x,n)$ is measurable with
respect to the $\sigma $-algebra $\mathcal G_{n}^{n+k+1}$, where
\begin{align}
  \label{MBNG-eq:defGnnk}
  \mathcal{G}_{n}^{m} := 
  \sigma\big(\omega(y,i), \widetilde\omega{(y,i)}  :
    y \in \ZZ^d, n \leq i < m \big), \quad n < m.
\end{align}
For $k\in\{0,\dots,\infty\}$, we define $M_k(x,n)\subseteq U(x,n)$ to be the set of
sites which maximise $\ell_k$ over $U(x,n)$, i.e.\
\begin{equation*}
  M_k(x,n) := \Big\{y\in U(x,n):\ell_k(y)=\max_{z\in U(x,n)}
  \ell_k(z)\Big\},
\end{equation*}
and for convenience we set $M_{-1}(x,n)=U(x,n)$.
Observe that we have
\begin{align*}
  M_0(x,n)&=\{y\in U(x,n): y\text{ is open}\},\\
  M_\infty(x,n)&= U(x,n)\cap \mathcal C, \\
  M_k(x,n)&\supseteq M_{k+1}(x,n), \qquad k\ge -1.
\end{align*}
Let $m_k(x,n)\in M_k(x,n)$ be the element of $M_k(x,n)$ that appears as the
first in the permutation $\widetilde\omega{(x,n)}$.

\begin{figure} 
 \centering
 \begin{tikzpicture}[xscale=0.62,yscale=0.62]
   \draw (0,-0.5) node {\footnotesize $(x,n)$};
   \draw (2.1,0.65) node {\footnotesize $\widetilde \omega_{(x,n)}(1)$};
   \draw (-2.05,0.65) node {\footnotesize $\widetilde \omega_{(x,n)}(2)$};
   \filldraw[fill=black] (0,0) circle (3pt)
                         (-1,1) circle (3pt)
                         (1,1) circle (3pt)
                         (2,2) circle (3pt)
                         (-2,2) circle (3pt)
                         (-3,3) circle (3pt)
                         (-1,3) circle (3pt);
   \filldraw[fill=white] (0,2) circle (3pt)
                         (3,3) circle (3pt)
                         (1,3) circle (3pt);

   \draw[-stealth,thick] (0.1,0.1) -- (0.9,0.9);
   \draw[-stealth,thick] (-0.9,1.1) -- (-0.1,1.9);
   \draw[-stealth,thick] (0.9,1.1) -- (0.1,1.9);
   \draw[-stealth,thick] (-2.1,2.1) -- (-2.9,2.9);
   \draw[-stealth,thick] (0.1,2.1) -- (0.9,2.9);
   \draw[-stealth,thick] (2.1,2.1) -- (2.9,2.9);

   \draw[-stealth,thick] (3.1,3.1) -- (3.9,3.9);
   \draw[-stealth,thick] (-2.9,3.1) -- (-2.1,3.9);
   \draw[-stealth,thick] (-0.9,3.1) -- (-0.1,3.9);
   \draw[-stealth,thick] (0.9,3.1) -- (0.1,3.9);

   \draw[-stealth,thick,dotted] (-0.1,0.1) -- (-0.9,0.9);
   \draw[-stealth,thick,dotted] (-1.1,1.1) -- (-1.9,1.9);
   \draw[-stealth,thick,dotted] (-1.9,2.1) -- (-1.1,2.9);
   \draw[-stealth,thick,dotted] (-0.1,2.1) -- (-0.9,2.9);
   \draw[-stealth,thick,dotted] (1.1,1.1) -- (1.9,1.9);
   \draw[-stealth,thick,dotted] (1.9,2.1) -- (1.1,2.9);
   \draw[-stealth,thick,dotted] (-3.1,3.1) -- (-3.9,3.9);
   \draw[-stealth,thick,dotted] (-1.1,3.1) -- (-1.9,3.9);
   \draw[-stealth,thick,dotted] (1.1,3.1) -- (1.9,3.9);
   \draw[-stealth,thick,dotted] (2.9,3.1) -- (2.1,3.9);

   \draw[thick] (4,0) -- (5,1);
   \filldraw[fill=black] (4,0) circle (1.5pt) (5,1) circle (1.5pt);
   \draw (4.5,-0.4) node {\footnotesize $k=1$};


   \draw[thick] (7,0) -- (8,1) -- (7,2);
   \filldraw[fill=black] (7,0) circle (1.5pt) (8,1) circle (1.5pt) (7,2) circle (1.5pt);
   \draw (7.5,-0.4) node {\footnotesize $k=2$};

   \draw[thick] (9,0) -- (10,1) -- (11,2) -- (12,3) ;
   \filldraw[fill=black] (9,0) circle (1.5pt) (10,1) circle (1.5pt) (11,2) circle
   (1.5pt) (12,3) circle (1.5pt);
   \draw (10.5,-0.4) node {\footnotesize $k=3$};

   \draw[thick] (16,0) -- (15,1) -- (14,2) -- (13,3) -- (14,4) ;
   \filldraw[fill=black] (16,0) circle (1.5pt) (15,1) circle (1.5pt) (14,2) circle
   (1.5pt) (13,3) circle (1.5pt) (14,4) circle (1.5pt);
   \draw (14,-0.4) node {\footnotesize $k=4$};
 \end{tikzpicture}
 \caption{The paths $\gamma_k^{(x,n)}$ from \eqref{MBNG-e:gamma} based on
   $\omega$'s and $\widetilde\omega$'s. Black and white circles represent
   open sites, i.e.\ $\omega(\text{site})=1$, and closed sites, i.e.\
   $\omega(\text{site})=0$, respectively. Solid arrows from a site point
   to $\widetilde\omega_{(\text{site})}(1)$ and dotted to
   $\widetilde\omega_{(\text{site})}(2)$. On the right the sequence of
   paths $\gamma_k^{(x,n)}(\cdot)$ for $k=1,2,3,4$ is shown.
   For sake of pictorial
   clarity, we used here $U(x,n)=\{(x+1,n+1),(x-1,n+1)\}$ here instead of \eqref{MBNG-Udef}.}
 \label{MBNG-fig:gammapaths}
\end{figure}
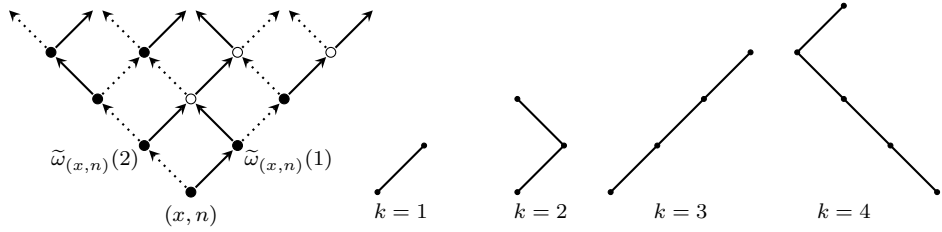

Given $(x,n)$, $k$, $\omega $ and $\widetilde\omega$, we define a path
$\gamma_k = \gamma_k^{(x,n)}$ of length $k$ via
\begin{equation}
  \label{MBNG-e:gamma}
  \gamma_k(0) = (x,n), \qquad 
    \gamma_k(j+1) = m_{k-j-2}(\gamma_k(j)) \;\;\text{for } j=0,\dots,k-1.
\end{equation}
In words, at every step, $\gamma_k$ checks the neighbours of its
present position and picks randomly (using the random permutation
$\widetilde \omega $) one of those where it can go further on open
sites, but inspecting only the state of sites in the time-layers
$\{n,\dots,n+k-1\}$. Consequently, the construction of
$\gamma_k^{(x,n)}$ is measurable with respect to the $\sigma $-algebra
$\mathcal G_n^{n+k}$ from \eqref{MBNG-eq:defGnnk}. See
Figure~\ref{MBNG-fig:gammapaths} for an illustration. 
Intuitively, $\gamma_k^{(x,n)}$ would be the trajectory of $X$ starting from the 
space-time point $(x,n)$ if we replaced in 
\eqref{MBNG-eq:defXdynamics} the condition that $X$ can only walk on $\cC$ 
by the requirement that the first $k$ steps must begin on open sites. 

It is not hard to check that these paths $\gamma_k^{(x,n)}$ have the following properties 
(see \cite[Lemma~2.1 and Remark~2.2]{MBNG-BCDG13} for details): 
Given $\omega$, $(x,n)\in \cC$ and $\widetilde\omega$, 
\begin{itemize}
\item[(a)] (steps begin on open sites) $\omega(\gamma_k(m))=1$ for all $0\le m < k$.
\item[(b)] (stability in $k$)
  If the end point of $\gamma_k$ is open,
  i.e.~$\omega (\gamma_k(k))=1$, then the path
  $\gamma_{k+1}$ restricted to the first $k$ steps equals
  $\gamma_k$.
\item[(c)] (fixation on $\mathcal C$)
  Assume that $\gamma_k(j)\in \mathcal C$ for some $k\ge 0$, $j\le k$.
  Then, $\gamma_m(j)=\gamma_k(j)$ for all $m>k$.
\item[(d)] (exploration of finite branches) If
  $\gamma_k(k-1)\in \mathcal C $ and $\gamma_k(k)\notin\mathcal  C$ for
  some $k$, then $\gamma_j(k)=\gamma_k(k)$ for all
  $k\le j \le k+\ell(\gamma_k(k))+1$ and
  $\gamma_{k+\ell(\gamma_k(k))+2}(k)\neq \gamma_k(k)$.
\end{itemize}
By (c), $\gamma_\infty^{(x,n)}(j) = \lim_{k\to\infty} \gamma_k(j)$ exists a.s.\ 
(since holes in the cluster are a.s.\ finite).
Furthermore, for fixed $\omega$ and $(x,n)\in \cC$ (but thinking of $\widetilde\omega$ as random), 
the law of $(\gamma_\infty^{(x,n)}(j))_{j\ge 0}$ is the same as the law of the random
walk $(X_j,n+j)_{j\ge 0}$ on $\mathcal C$ started from $(x,n)$.
Thus we can and shall couple the random walk $(X_k,k)$ started from
  $(\boldsymbol 0,0)$ with the random variables
  $\omega, \widetilde\omega$ by setting
  \begin{equation}
    \label{MBNG-eq:walkpathcoupled}
    (X_k,k)= \gamma_\infty^{(\boldsymbol 0,0)}(k)
    =\lim_{j\to\infty} \gamma_j^{(\boldsymbol 0,0)}(k).
  \end{equation}

With these ingredients, we can define regeneration times as follows: 
Let 
\begin{align}
  \label{MBNG-eq:reg.1}
  T_0:= 0 \quad \text{and}\quad
  T_{j}:=\inf\left\{k> T_{j-1}: \xi(\gamma^{(\boldsymbol 0,0)}_k(k))=1\right\},\quad j\ge 1.
\end{align}
(Here and later we use the notation $\xi (y) := \xi_n(x)$ when
  $y=(x,n)\in\mathbb Z^d\times \mathbb Z$.)
At times $T_j$ the local construction of the path finds a `real ancestor'
of $(\boldsymbol 0,0)$ in the sense that for any $m>T_j$,
$\gamma^{(\boldsymbol 0,0)}_m(T_j)=\gamma^{(\boldsymbol 0,0)}_{T_j}(T_j)$, by
property~(c).
The increments between regeneration times are 
\begin{align}
\label{MBNG-deftauiYi}
\tau_i:= T_i-T_{i-1} \qquad \text{and} \qquad
Y_i:= X_{T_i}-X_{T_{i-1}}.
\end{align}
and 
we then indeed have that
\begin{align} 
  \label{MBNG-Ytauiiid}
  & \text{the sequence $\bigl((Y_i,\tau_i)\bigr)_{i \ge 1}$ is
  i.i.d.\ and $Y_{1}$ is symmetrically distributed}, \\[1ex]
  \label{MBNG-Ytautails}
  & \text{both $Y_1$ and $\tau_1$ have exponential tails}.
\end{align}
\index{regeneration construction!for directed random walk on oriented percolation}%

The intuition behind the regeneration property \eqref{MBNG-Ytauiiid} is the following: 
Assume that for some $k$, we have constructed the path $\gamma_k^{(\boldsymbol 0,0)}$ and observe 
that $\xi\big(\gamma_k^{(\boldsymbol 0,0)}(k)\big)=1$. Then we have obtained 
information about some $\omega(y,j)$ 
and $\widetilde\omega(y,j)$ for $j<k$, $y \in \ZZ^d$ and we know that the site $\gamma_k^{(\boldsymbol 0,0)}(k)$ 
in time-slice $k$ is connected to $+\infty$. The latter property depends only on $\omega(y,j)$ with $j \ge k$, $y \in \ZZ^d$ 
and the $\omega$'s in different time-slices are independent. By property (c), we have $(X_k,k) = \gamma_k^{(\boldsymbol 0,0)}(k)$.
Thus, concerning the future behaviour of $X$, we are then at time $k$ in the same situation as at time $0$: 
All we know (and need to know) is that $X$ sits on some site in $\cC$, and we can start afresh.

However, if we observe that $\xi\big(\gamma_k^{(\boldsymbol 0,0)}(k)\big)=0$, we are in a different situation: 
We then know that $\gamma_k^{(\boldsymbol 0,0)}(k)$ is the starting point of a finite (possibly empty) 
oriented percolation cluster. Then we must continue the local construction until it has explored the  
`reason why $\xi\big(\gamma_k^{(\boldsymbol 0,0)}(k)\big)=0$', which depends on finitely many sites
(cf property~(d) above). 

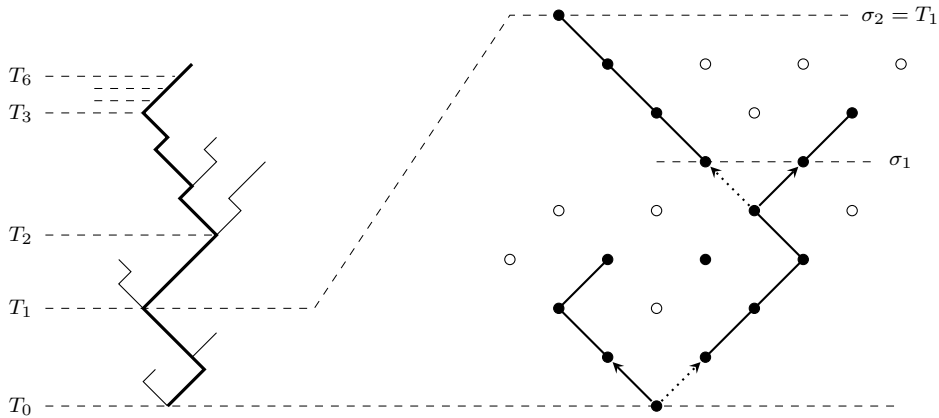
\begin{figure} 
  \centering
  \begin{tikzpicture}[xscale=0.65,yscale=0.65]

    \draw (-12,0) -- (-12.5,0.5)--(-12.25,0.75);
    \draw (-11.5,1) -- (-11,1.5);
    \draw (-12.5,2) -- (-13,2.5) -- (-12.75,2.75) -- (-13,3);
    \draw (-11,3.5) -- (-10.5,4) -- (-10.75,4.25) -- (-10,5);
    \draw (-11.5,4.5) -- (-11,5)--(-11.25,5.25)--(-11,5.5) ;

    \draw[very thick]  (-12,0)-- (-11.25,0.75) --
    (-12.5,2)--(-12,2.5)--(-11,3.5) --(-11.75,4.25) --
    (-11.5,4.5) -- (-12.25,5.25) -- (-12,5.5) -- (-12.5,6) -- (-11.5,7)  ;


    \draw[dashed, very thin] (-14.5,0) -- (3,0);   \draw (-15,0) node {\footnotesize $T_0$};
    \draw[dashed, very thin] (-14.5,2) -- (-9,2) -- (-5,8) -- (2.0,8);   \draw (-15,2) node {\footnotesize $T_1$};

    \draw[dashed, very thin] (-14.5,3.5) -- (-11.01,3.5);   \draw (-15,3.5) node {\footnotesize $T_2$};
    \draw[dashed, very thin] (-14.5,6) -- (-12.55,6);  \draw (-15,6) node {\footnotesize $T_3$};
    \draw[dashed, very thin] (-13.5,6.25) -- (-12.20,6.25);
    \draw[dashed, very thin] (-13.5,6.5) -- (-12.10,6.5);
    \draw[dashed, very thin] (-14.5,6.75) -- (-11.85,6.75);  \draw (-15,6.75) node {\footnotesize $T_6$};


    \draw[dashed, very thin] (-2,5) -- (2.5,5);  \draw (3.0,5) node {\footnotesize $\sigma_1$};
    \draw[anchor=west] (2,8) node {\footnotesize $\sigma_2=T_1$};

    \filldraw[fill=black]
    (-2,0) circle (3pt)
    (-1,1) circle (3pt)
    (-3,1) circle (3pt)
    (-4,2) circle (3pt)
    (0,2) circle (3pt)
    (-3,3) circle (3pt)
    (1,3) circle (3pt)
    (-1,3) circle (3pt)
    (1,5) circle (3pt)
    (0,4) circle (3pt)
    (-1,5) circle (3pt)
    (2,6) circle (3pt)
    (-2,6) circle (3pt)
    (-3,7) circle (3pt)
    (-4,8) circle (3pt);
    

    \filldraw[fill=white]
    (-2,2) circle (3pt)
    (-5,3) circle (3pt)
    (-4,4) circle (3pt)
    (-2,4) circle (3pt)
    (2,4) circle (3pt)
    (0,6) circle (3pt)
    (3,7) circle (3pt)
    (1,7) circle (3pt)
    (-1,7) circle (3pt);


    \draw[thick] (-3,1) -- (-4,2) -- (-3,3);
    \draw[thick] (-1,1) -- (1,3)-- (0,4);
    \draw[thick] (1,5) -- (2,6);
    \draw[thick] (-1,5) -- (-4,8);

    \draw[-stealth,thick] (-2.1,0.1) -- (-2.9,0.9);
    \draw[-stealth,thick,dotted] (-1.9,0.1) -- (-1.1,0.9);

    \draw[-stealth,thick] (0.1,4.1) -- (0.9,4.9);
    \draw[-stealth,thick,dotted] (-0.1,4.1) -- (-0.9,4.9);




  \end{tikzpicture}
  \caption{`Discovering' of the trajectory of $X$ between the
    regeneration times $T_0$ and $T_6$ in case $U =\{-1,1\}$ is shown on
    the left-hand side of the figure. On the right-hand side we zoom into
    the evolution between $T_0$ and $T_1$. On the two `relevant sites'
    we show the realisation of $\widetilde \omega$'s using the same
    conventions as in Figure~\ref{MBNG-fig:gammapaths} (in particular,
      again $U(x,n)=\{(x+1,n+1),(x-1,n+1)\}$).}
  \label{MBNG-fig:expl-branch}
\end{figure}

See Figure~\ref{MBNG-fig:expl-branch} for an illustration: In this example, the local construction enters 
a finite cluster at time $k=\sigma_1$ and explores this, regeneration occurs then at time $T_1=\sigma_2$ 
when the exploration is completed. The full details are in \cite[Lemma~2.5]{MBNG-BCDG13}.
\smallskip

To obtain \eqref{MBNG-Ytautails}, one uses the fact that the height of a finite cluster
in supercritical oriented percolation has exponential tails, see
\cite{MBNG-Durrett:84} and \cite[Lemma~A.1]{MBNG-BCDG13}.
The distributional symmetry of $Y_1$ follows from the symmetry of $U(x,n)$ in \eqref{MBNG-Udef}.

\medskip

Given \eqref{MBNG-Ytauiiid} and \eqref{MBNG-Ytautails}, 
the law of large numbers \eqref{MBNG-eq:LLN} and the annealed CLT \eqref{MBNG-eq:annealedCLT} 
follow straightforwardly by re-writing $X_n$ as a sum along regeneration times plus an 
asymptotically negligible remainder. 
The quenched CLT \eqref{MBNG-eq:quenchedCLT} requires some additional effort: 
Here, we used two copies $X$ and $X'$ of the walk on the same cluster $\cC$ 
to control the variance of $E_\omega\big[f\left(X_n/\sqrt{n} \, \right) \,\big]$, 
an approach inspired by \cite{MBNG-BolthausenSznitman:2002}. This, in turn, requires to 
enlarge the regeneration construction to incorporate simultaneous regenerations 
for both $X$ and $X'$. Studying two (or more) copies of the walk on $\cC$, especially when 
one stipulates that they coalesce as soon as they meet, is also 
very natural from the point of view of larger samples. 
In fact, this is exactly the device that is made use of in \cite{MBNG-BGS18} and
it also plays a key role in the proof of \eqref{MBNG-eq:tailbounds} below. 
We will however not spell out the details here and instead refer to \cite{MBNG-BCDG13, MBNG-BGS18}.

\subsection{Extensions}

\subsubsection{Contact process with fluctuating population sizes}

Let $K(x,n)$, $(x,n) \in \ZZ^d \times \ZZ$ be possibly correlated $\NN$-valued
  random variables, independent of the $\omega$'s.  We define the {discrete time contact process
    with fluctuating population size},
  $\widehat \eta := (\widehat\eta_n)_{n \in \ZZ}$, by
  \begin{align}
    \label{eq:8}
    \widehat \eta_n(x) := \eta_n(x) K(x,n),
  \end{align}
  with $\eta_n(x)$ from \eqref{MBNG-eq:5} and its time reversal $\widehat \xi := ( \widehat \xi_n)_{n \in \ZZ}$ by
    $\widehat \xi_n(x) := \xi_n(x) K(x,n).$
    
  One can interpret $K(x,n)$ as a random `carrying capacity' of the site
  $(x,n)$: When $\eta_n(x)=1$, $K(x,n)$ individuals live at position $x$ in generation $n$,
  and each of them is independently assigned an ancestor from $\cA_{x,n}$ as
  in \eqref{MBNG-eq:defancestors}.
  \index{fluctuating population size}%
  
  Now conditioned on $\widehat\xi_0(\mathbf{0}) \ge 1$ the ancestral random
  walk 
  is defined by $X_0=\boldsymbol 0$ and \eqref{MBNG-eq:defXdynamics} is generalised to 
  \begin{equation}
    \label{MBNG-eq:11}
    \PP\bigl( X_{n+1}=y \mid X_{n}=x, \widehat \xi\, \bigr) =
    \begin{cases}
      \dfrac{\widehat \xi_{n+1}(y)} {\sum_{(y',n+1) \in U(x,n)}
        \widehat\xi_{n+1}(y')} \;\;
       &\text{if }(y,n+1) \in  U(x,n),\\[4mm]
      0&\text{otherwise.}
    \end{cases}
  \end{equation}
  \index{random walk!ancestral}%
  
  Analogues of Theorem~\ref{MBNG-thmrwopc} then hold under suitable
  assumptions on the random field
  $K=(K(x,n))_{x \in \ZZ^d, n\in \ZZ}.$ The case when $K$ is an
  i.i.d.\ field is discussed in \cite[Remark~1.6]{MBNG-BCDG13}.
  K.~Miller \cite{MBNG-M16} generalises this considerably by assuming
  instead certain mixing conditions: A law of large numbers analogous
  to \eqref{MBNG-eq:LLN}, with possibly non-zero speed, holds if $K$
  is $\phi$-mixing in time with coefficients
  $\phi_n \in O(n^{-1-\delta})$ for some $\delta>0$, an annealed CLT
  analogous to \eqref{MBNG-eq:annealedCLT} holds if
  $\phi_n \in O(n^{-2-\delta})$; a quenched CLT analogous to
  \eqref{MBNG-eq:annealedCLT} holds if $K$ is exponentially mixing in
  space and time. Note that in general, \eqref{MBNG-eq:11} describes
  a non-elliptic random walk in a non-Markovian (but mixing) environment.
  The key idea is again a `regeneration construction'
  where the i.i.d.\ property in \eqref{MBNG-Ytauiiid} is now replaced
  by a sufficiently strong mixing property. We refer to \cite{MBNG-M16} and \cite{MBNG-M17}
  for details. 
 

\medskip

\subsubsection{Brownian web limit in spatial dimension one}
\label{MBNG-subsect:BrownianWeb}

One can consider the ancestral lineages of all individuals in the stationary $\eta$ 
from \eqref{MBNG-eq:5} simultaneously. This gives rise to an infinite system of
random walks $X^{(x,n)} = (X^{(x,n)}_m)_{m\ge n}$ on the time-reversal $\xi$ of $\eta$, where for
each $(x,n) \in \cC$, the walk $X^{(x,n)}$ starts at time $n$ at position $x$, follows the
analogue of \eqref{MBNG-eq:defXdynamics}, and different walkers coalesce whenever they meet in the
same space-time site. By Theorem~\ref{MBNG-thmrwopc} and space-time stationarity, any $X^{(x,n)}$
converges to a Brownian motion under diffusive rescaling. As shown in \cite{MBNG-BGS18},
in spatial dimension $d=1$, the collection of all these paths converges after diffusive
rescaling as in Theorem~\ref{MBNG-thmrwopc} in distribution
to the Brownian web. Informally, this limit object describes an infinite system of coalescing
Brownian motions starting from all space-time points in $\RR \times \RR$.
One may then apply our convergence result to investigate the
behaviour of interfaces in the discrete time contact process analogously to \cite[Theorem 7.6
and Remark 7.7]{MBNG-NewmanRavishankarSun2005}, as observed in \cite[p.~1051]{MBNG-BGS18}. 
\index{Brownian web}\index{interface}%
We refer also to the article `Interfaces in spatial population dynamics'
by Marcel Ortgiese in this volume, 
which studies spatial population models (in continuous
space) in $d=1$, with a particular focus on interfaces.  These models
are `continuum analogues' of the voter model, and the interfaces are
stochastic processes in dynamic environments. Dualities and their
genealogical interpretations play an important role there as well.

An important ingredient in the proof is a quantitative strengthening 
of \eqref{MBNG-multitypd12} from Remark~\ref{MBNG-rem:multitpye}: 
\begin{align}
  \label{MBNG-eq:tailbounds} 
  \PP\left(T^{(z_1, z_2)}_{\mathrm{meet}} > n \, \Big| \, 
  \xi_0(z_1)=\xi_0(z_2)=1 \right) \le \text{const.}\times \frac{|z_1-z_2|}{\sqrt{n}}  
  \quad \text{for } z_1, z_2 \in \ZZ, \, n \in \NN, 
\end{align}
\index{meeting time!for two random walks on the oriented percolation cluster}%
where $T^{(z_1, z_2)}_{\mathrm{meet}}$ is the number of steps until two walks on the same realisation 
of $\xi$ which start at time $0$ from $z_1$ and $z_2$, respectively, meet for the first time.
Note that \eqref{MBNG-eq:tailbounds} is the asymptotically correct form of the 
decay for simple random walks in $d=1$.
For more information, we refer to \cite{MBNG-S17}.

The results in \cite{MBNG-BGS18} can again be interpreted as
an averaging statement about the percolation cluster: apart from a
change of variance, it behaves as the full lattice (for which convergence
to the Brownian web was proved in \cite{MBNG-NewmanRavishankarSun2005}), i.e. the effect of
the `holes' in the cluster vanishes on a large scale.
For a thorough discussion of the Brownian web, including historical comments and references, see 
the overview article \cite{MBNG-SSS17}.
Note that there is no analogous object in spatial dimension $d \ge 2$ because there,
independent Brownian motions never meet.

\section{Ancestral lineages for logistic branching random walks}
\label{MBNG-sect:lbrw}

We consider a system of discrete-time branching random walks with logistic regulation: 
Let $\eta_n(x)$ be the number of individuals at position $x \in \ZZ^d$ in generation $n \in \ZZ$. 
Given the configuration $\eta_n$ at time $n$,
for $x \in \ZZ^d$, each individual at $x$ has a Poisson-distributed number of offspring
with mean
\begin{align}
  \label{MBNG-eq:logisticmean}
  \big(m - \sum_z \lambda_{z-x} \eta_n(z)\big)^+
\end{align}
\index{logistic branching random walk}%
and each child moves to $y$ with probability $p_{y-x}$, independently
for different parental individuals and for different children.  Here,
$p_{xy} = p_{y-x}$ is a symmetric, aperiodic finite range random walk
kernel on $\ZZ^d$, $m > 1$, $\lambda_z \ge 0$, $z \in \ZZ^d$ is
symmetric with finite range and $\lambda_0>0$.  These children then
form the next generation, $\eta_{n+1}$.  \eqref{MBNG-eq:logisticmean}
has a natural interpretation as local competition: each individual at
$z$ reduces the average reproductive success of a focal individual at
$x$ by $\lambda_{z-x}$. In particular, this introduces local
density-dependent feedback in the model: The offspring distribution
is supercritical when there are few neighbours and subcritical when there
are many neighbours. 
Note that by properties of the Poisson distribution $(\eta_n)$ is in fact a 
probabilistic cellular automaton: Given $\eta_n$,
\begin{align}
  \label{MBNG-eq:PCA}
  \eta_{n+1}(y) \sim \mathrm{Poisson}\left(
  {\sum\limits_{x \in \ZZ^d} } \big(m - \sum_{z \in \ZZ^d} \lambda_{z-x} \eta_n(z)\big)^+ p_{y-x}\right) ,
\end{align}
independently for different $y \in \ZZ^d$.
\index{probabilistic cellular automaton}%
\index{branching random walk}%
\smallskip

\begin{remark} 1.\ For the choice $\lambda \equiv 0$, the system
  $(\eta_n)_n$ is a `classical' branching random walk, in which
  different individuals behave completely independently.
  This is a classical topic with a lot of recent progress, 
  see in particular the article `Branching random walks in random environment' by
  Wolfgang K\"{o}nig in this volume.
  In \cite{MBNG-GKS:13} and \cite{MBNG-GKS:15},
  moment asymptotics for the number of particles in a branching random
  walk in random environment are derived. Note that the first moments
  correspond to the well-investigated solutions of the parabolic
  Anderson model.\smallskip %

  \noindent 2.\ Conditioning on $\eta_n(\cdot) \equiv N$ in
  \eqref{MBNG-eq:PCA} for some $N \in \NN$ 
  and considering types and/or ancestral
  relationships, as we will do below, yields a 
  version of the stepping stone model.  \smallskip

  \noindent 3.\ The form of the competition kernel and the Poisson offspring
  law in \eqref{MBNG-eq:logisticmean}--\eqref{MBNG-eq:PCA} are
  prototypical (and convenient for the proofs) but can be replaced by
  more general choices, see the discussion in
  \cite[Remark~5~(ii)]{MBNG-BD07} and \cite[Section~5]{MBNG-BCD16}.
\end{remark}

\begin{theorem}[Survival and complete convergence, {\cite[Theorem~1 and Corollary~4]{MBNG-BD07}}]
  \label{MBNG-thm:logbrw}
  Assume $m \in (1,3)$. There exist $\varepsilon_0, \varepsilon_1 >0$ such that 
  for all choices $0< \lambda_0 \le \varepsilon_0$ and $0 \le \lambda_z \le \varepsilon_1 \lambda_0$ 
  for $z \neq 0$, the system $(\eta_n)$ survives for all time locally (and hence also globally)
  with positive probability for any non-trivial initial condition $\eta_0$.
  %
  Given survival (either local or global), $\eta_n$ converges as $n\to\infty$ in distribution to 
  its unique non-trivial equilibrium.
\end{theorem}
\index{complete convergence}%
\index{survival}%

\begin{figure}
  \begin{center}
    \includegraphics[width=6.4cm]{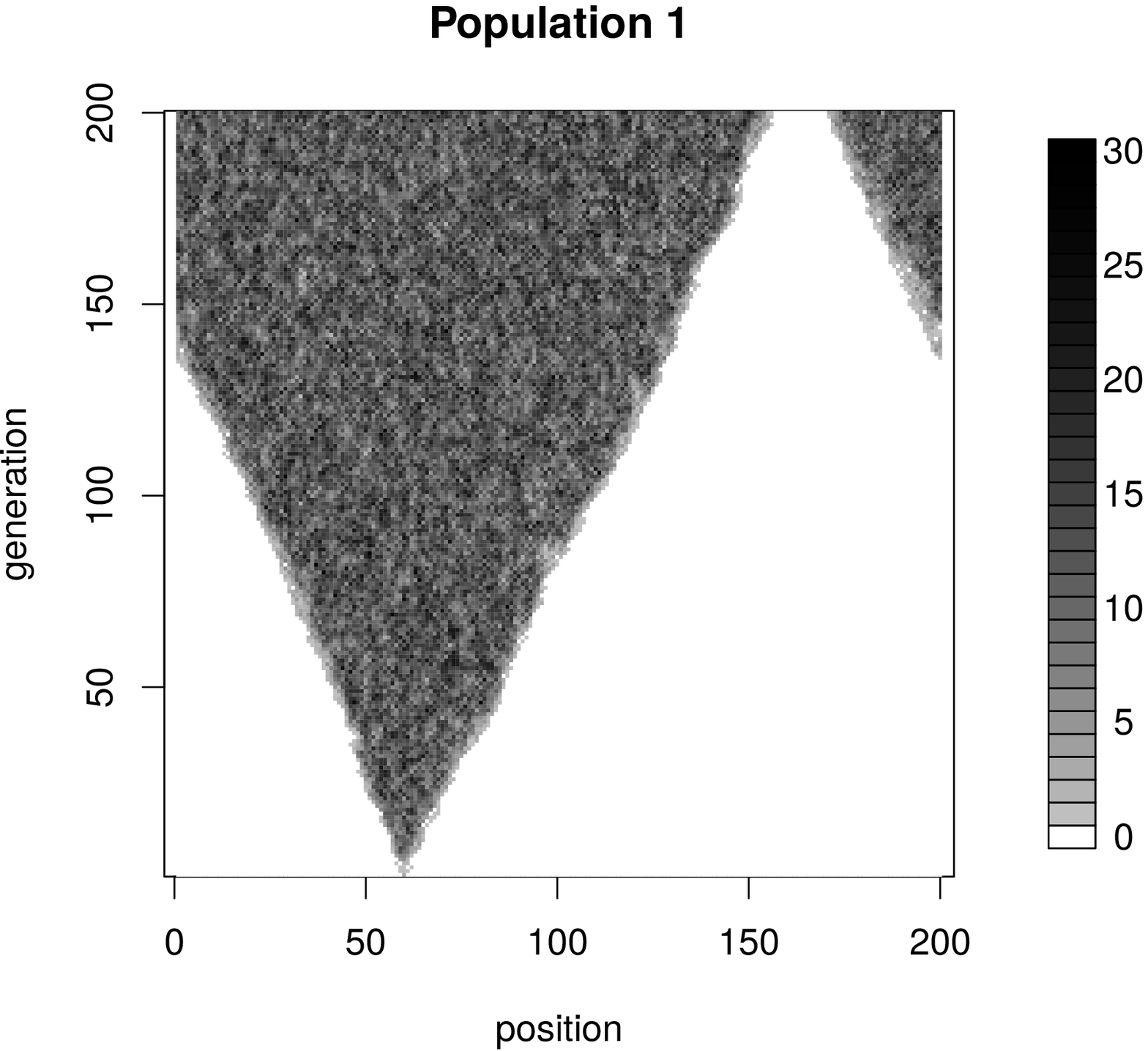}%
    \includegraphics[width=6.4cm]{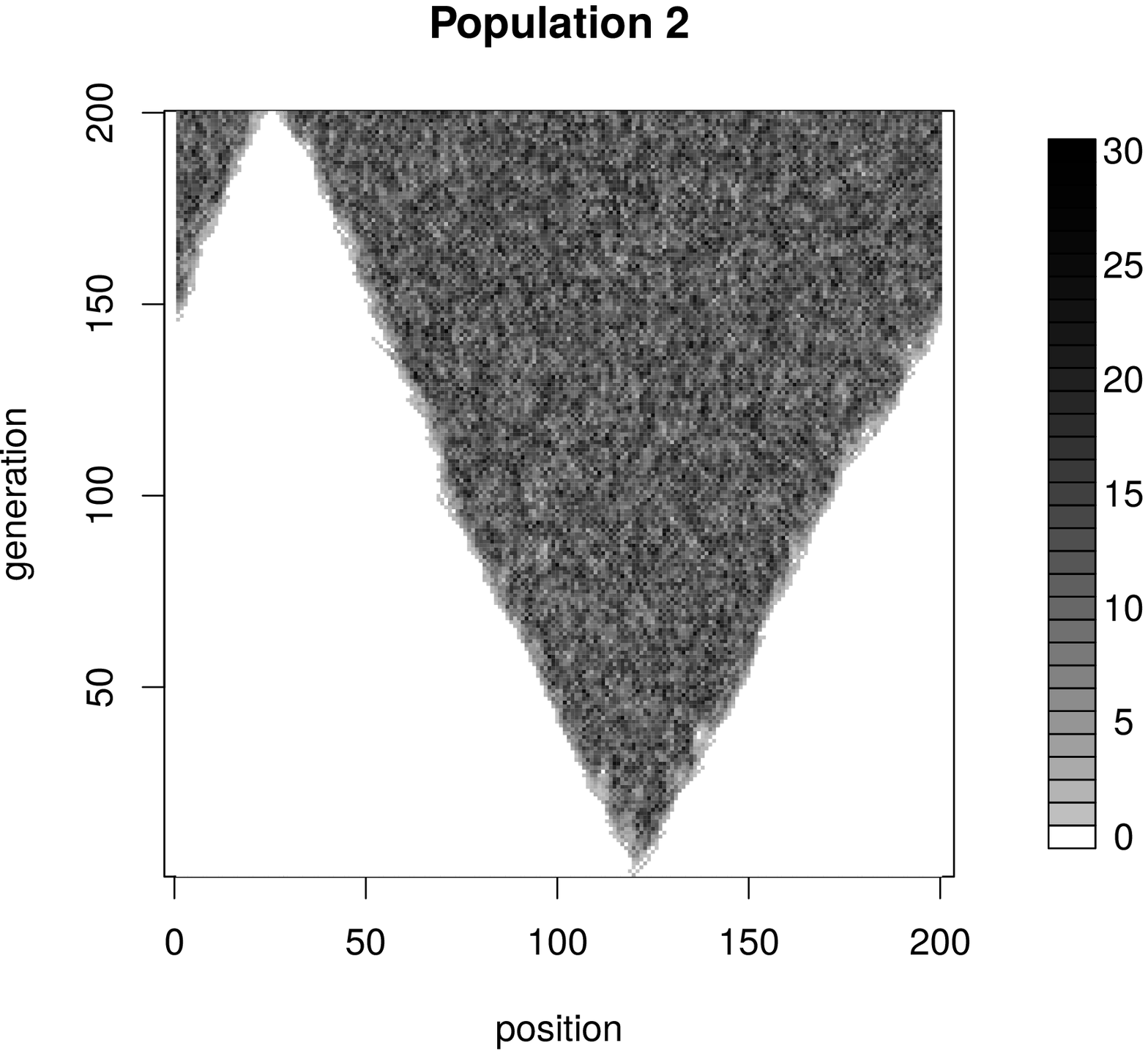}\\[-0.5cm]
    \includegraphics[width=6.4cm]{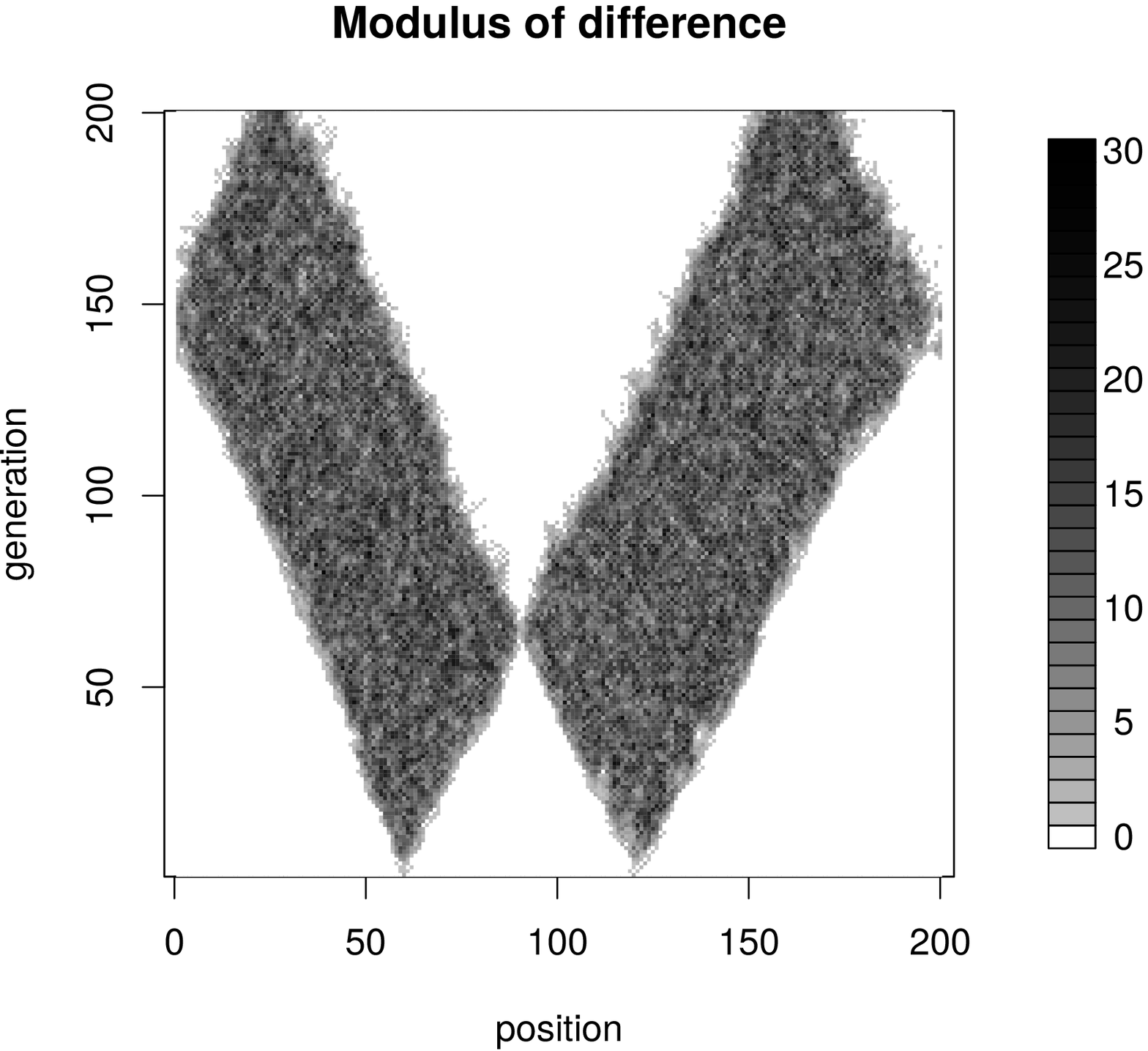}
  \end{center}
  \caption{Starting from any two initial conditions $\eta_0$, $\eta_0'$, 
    copies $(\eta_n)$, $(\eta'_n)$ can be coupled such that 
    if both survive (here, $m=1.5$, $p=(1/3,1/3,1/3)$, $\lambda=(0.01,0.02,0.01)$,
    $\eta_0=\delta_{60}$, $\eta_0'=\delta_{120}$ and space is $\{1,2,\dots,200\}$
    with periodic boundary conditions).
    The picture at the bottom shows $|\eta_n(x) - \eta_n'(x)|$, note the 
    growing region in the middle where $\eta_n(x)=\eta_n'(x)>0$.}
  \label{MBNG-fig:coupling}
\end{figure}

We will not prove Theorem~\ref{MBNG-thm:logbrw} here but point out that a crucial ingredient in
the proof is a strong coupling property of the system $(\eta_n)$: Starting from any two initial
conditions $\eta_0$, $\eta_0'$, 
\begin{equation}
\label{MBNG-lbrwcouplingclaim}
\parbox{0.8\textwidth}{copies $(\eta_n)$, $(\eta'_n)$ can be
coupled such that if both survive, \\[0.5ex]
$\eta_n(x)=\eta'_n(x)$ in a
space-time cone.}
\end{equation}
This allows to compare the system to supercritical oriented
percolation on suitably coarse-grained space-time scales, see
\cite[Section~5]{MBNG-BD07} for details and see
Figure~\ref{MBNG-fig:coupling} for a simulation.
\index{coupling}%

\smallskip

\begin{remark} 1.\ 
  The restriction to $m<3$ in Theorem~\ref{MBNG-thm:logbrw} is
  `inherited' from the logistic iteration $w_{n+1} = m w_n(1-w_n)$
  because in this parameter regime, it has a unique attracting fixed
  point. Note that literally, a `deterministic space-less' analogue of
  \eqref{MBNG-eq:PCA} would read
  $\widetilde{w}_{n+1} = m \widetilde{w}_n(1- \tilde{\lambda}
  \widetilde{w}_n)$ with $\tilde{\lambda} = \sum_{z} \lambda_z$, the
  rescaling $\widetilde{w}_n = (m/\tilde{\lambda})w_n$ brings this to
  the `standard form' just mentioned.
  
  Survival can be proved also for $m \in [3,4)$ with similar
  arguments, but convergence cannot.  For $m < 1$ (and for $m=1$ in
  $d \leq 2$) one can easily see, using domination by subcritical
  branching random walks, that $(\eta_n)_n$ will die out locally when
  starting from any initial condition $\eta_0$ with
  $\sup_{x \in \ZZ^d} \EE[\eta_0(x)] < \infty$.
  \smallskip

  \noindent 2.\ In \cite{MBNG-GrevenSturmWinterZaehle:19+}, multi-type
  continuous mass branching populations with competitive interactions
  are studied, the logistic branching random walks we described in
  Section~\ref{MBNG-sect:lbrw} are a close relative of such systems in
  the single-type case (or in the multi-type case with completely
  symmetric parameters). Furthermore, by using space $\RR^d$ as a
  `trait space,' the measure-valued processes studied in
  \cite{MBNG-BSW18} can be seen as a suitable scaling limit of
  (relatives of) logistic branching random walks, see
  \cite[Remark~5]{MBNG-BSW18}.  Many challenging questions about the
  long-time behaviour of such continuous-mass interacting multi-type
  systems remain open. It is conceivable that the regeneration
  constructions for ancestral lineages we investigated in 
  \cite{MBNG-BCDG13} and \cite{MBNG-BCD16} might be adaptable to this
  context, and that this could enrich the pertinent `tool box.'
\end{remark}

\subsection{Dynamics of an ancestral lineage}

By Theorem~\ref{MBNG-thm:logbrw}, for suitable choices of the
parameters, the system~\eqref{MBNG-eq:PCA} has a unique non-trivial
equilibrium. We denote by
$\eta^{\rm stat} = (\eta^{\rm stat}_n(x))_{n \in \ZZ, x \in \ZZ^d}$
the corresponding stationary process and -- implicitly in our notation
-- `enrich' it suitably to allow bookkeeping of genealogical
relationships, as described at the beginning of
Section~\ref{MBNG-sect:lbrw}.  Consider the stationary
$\eta^{\rm stat}$ conditional on $\eta^{\rm stat}_0(\mathbf{0})>0$ and
sample an individual (uniformly) from the space-time origin
$(\mathbf{0},0)$, let $X_n$ be the spatial position of her ancestor
$n$ generations ago. Then
  \begin{align}
    \label{MBNG-eq:lbrwXtransprob}
    \PP\big( & X_{n+1}=y \, \big| \, X_n=x, \eta^{\rm stat} \big)
         = \frac{p_{x-y} \eta^{\rm stat}_{-n-1}(y) \big( m 
    - {\textstyle \sum_z \lambda_{z-y} \eta^{\rm stat}_{-n-1}(z) } \big)^+}{%
    \sum_{y'} p_{x-y'} \eta^{\rm stat}_{-n-1}(y') \big( m 
    - {\textstyle \sum_z \lambda_{z-y'} \eta^{\rm stat}_{-n-1}(z) } \big)^+} ,
  \end{align}
see \cite[(4.10--4.11)]{MBNG-BCD16}.
\index{logistic branching random walk!ancestral lineage in}%

Thus $(X_n)_n$ is a random walk in a -- relatively complicated -- 
random environment. Note that the forwards in time direction for the 
walk corresponds to backwards in time for $\eta^{\rm stat}$. Again it turns out that $X$ behaves 
like ordinary random walk when viewed over large enough space-time scales, as the 
following theorem shows.

\begin{theorem}[Law of large numbers and (averaged) central limit theorem, {\cite[Theorem~4.3]{MBNG-BCD16}}]
  \label{MBNG-thm:lbrw}
  Assume $m \in (1,3)$. There exist $\varepsilon_0, \varepsilon_1 >0$ such that 
  for all choices $0< \lambda_0 \le \varepsilon_0$ and $0 \le \lambda_z \le \varepsilon_1 \lambda_0$ 
  for $z \neq 0$, we have 
  \begin{equation}
    \label{MBNG-eq:lbrwLLN}
    \PP\Big( \frac1n X_n \to \mathbf{0} \,  \Big| \, \eta_0^{\rm stat}(\mathbf{0}) \neq 0 \Big) = 1 
  \end{equation}
  and 
  \begin{equation}
    \label{MBNG-eq:lbrwaCLT}
    \EE\Big[ f\big(\tfrac1{\sqrt{n}} X_n\big) 
    \, \Big| \, \eta_0^{\rm stat}(\mathbf{0}) \neq 0 \Big] 
    \mathop{\longrightarrow}_{n\to\infty} \EE\big[f(Z)\big]
  \end{equation}
  for $f \in C_b(\RR^d)$, where $Z$ is a (non-degenerate)
  $d$-dimensional normal random variable. A functional version of
  \eqref{MBNG-eq:lbrwaCLT} holds as well.
\end{theorem}
\index{law of large numbers}
\index{central limit theorem}

Note that \eqref{MBNG-eq:lbrwaCLT} is an averaged limit result. In ongoing
work with Andrej Depperschmidt and Timo Schl\"{u}ter, we are proving the
corresponding `quenched' limit theorem.

The proof of Theorem~\ref{MBNG-thm:lbrw} employs again a regeneration construction 
and a decomposition as in \eqref{MBNG-deftauiYi}. We will only sketch the main ideas 
below, referring the reader to \cite{MBNG-BCD16} for details.

Given $\eta^{\rm stat}$, $X$ is a time-inhomogeneous 
Markov chain; given also $X_n=x$ its transition probabilities in the $n+1$-th step 
depend only on $\eta^{\rm stat}_{-n-1}(x)$ in some finite window around $x$. 
We see from \eqref{MBNG-eq:lbrwXtransprob} that these transition 
probabilities are close to the fixed reference law $(p_x)_x$ if $X_n$ 
is in a region where the relative variation of $\eta_{-n-1}(X_n)$ is small. 

Thus, we define `good' space-time blocks in $\eta^{\rm stat}$ on suitable length 
scales $L_{\rm space} \gg 1$ and $L_{\rm time} \gg 1$, that is a finite set 
$\mathcal{G}$ of local configurations on $\{1,2,\dots,L_{\rm space}\}^d \times \{1,2,\dots,L_{\rm time}\}$ 
with the properties that 
\begin{itemize}
\item[(a)] $\eta^{\rm stat}$ has small relative variations inside a good block, 
\item[(b)] if the block with (coarse grained) index $(\tilde x,\tilde n) \in \ZZ^d \times \ZZ$ 
  is good, this will with high probability also be the case for its 
  `temporal successors' with indices $(\tilde x-1,\tilde n+1), (\tilde x,\tilde n+1),(\tilde x+1,\tilde n)$,
\item[(c)] if we consider two copies $\eta$ and $\eta'$ of the system \eqref{MBNG-eq:PCA} 
  with the property that in both the block with (coarse grained) index $(\tilde x,\tilde n)$ 
  is good, then with high probability the coupling discussed in \eqref{MBNG-lbrwcouplingclaim} 
  will make $\eta$ and $\eta'$ identical on the block with index $(\tilde x,\tilde n+1)$.
\end{itemize} 
(a) allows to control the walk $X$ whenever it moves through good blocks; 
(b) allows to compare the good blocks to supercritical oriented percolation (on the coarse-grained 
scale); (c) allows to `localise' information about the space-time configuration $\eta^{\rm stat}$ 
around good blocks, this is akin to the local construction from Section~\ref{MBNG-subsect:proofBCDG}. 
\smallskip

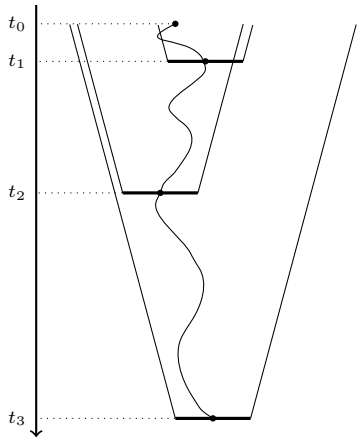
\begin{figure}
        \begin{tikzpicture}[scale=0.5]
        \draw [<-, thick] (-3,0) -- (-3,11.5);
        
        \filldraw[black] (0.7,11) circle (0.07cm); 
        \draw[dotted] (0.7,11) -- (-3,11); 
        \draw (-3.5,11) node {\scriptsize $t_0$};  
        \filldraw[black] (1.5,10) circle (0.07cm); 
        \draw[dotted] (1,10) -- (-3,10); 
        \draw (-3.5,10) node {\scriptsize $t_1$}; 
        \draw[very thick] (0.5,10) -- (2.5,10); 
        \draw (0.5,10) -- ++(105:1); 
        \draw (2.5,10) -- ++(75:1);
        
        \filldraw[black] (0.3,6.5) circle (0.07cm); 
        \draw[dotted] (0.3,6.5) -- (-3,6.5); 
        \draw (-3.5,6.5) node {\scriptsize $t_2$};  

        \draw[very thick] (-0.7,6.5) -- (1.3,6.5); 
        \draw (-0.7,6.5) -- ++(105:4.65);
        \draw (1.3,6.5) -- ++(75:4.65); 

        \filldraw[black] (1.7,0.5) circle (0.07cm); 
        \draw[dotted] (1.7,0.5) -- (-3,0.5); 
        \draw (-3.5,0.5) node {\scriptsize $t_3$};  

        \draw[very thick] (0.7,0.5) -- (2.7,0.5); 
        \draw (0.7,0.5) -- ++(105:10.85);
        \draw (2.7,0.5) -- ++(75:10.85); 

        \draw [rounded corners =3pt] 
        (0.7,11) .. controls (-0.5,10.4) and (0.9,10.7)..(1.5,10) 
        .. controls (1.5,9.5) and (0.3,9) .. (0.7,8.5) 
        .. controls (1.4,8) and (1,7.5).. (0.8,7.2) 
        .. controls (0.65,7) and (0.4,6.8) .. (0.3,6.5)
        .. controls (0,5.9) and (0.9,5.3) .. (1.2,4.7)
        .. controls (1.6,4.1) and (1.4,3.5) .. (1,2.9) 
        .. controls (0.6,2.2) and (0.9,1.5) .. (1.2,1) 
        .. controls (1.4,0.7) .. (1.7,0.5); 
      \end{tikzpicture}
      \caption{A schematic illustration of the regeneration
        construction for Theorem~\ref{MBNG-thm:lbrw}: The walk passes
        through a sequence of cones in an attempt to regenerate. Here, 
        regeneration at time $t_1$ fails because the path up to that 
        time is not covered by the corresponding cone and regeneration 
        at time $t_2$ fails because the corresponding cone does not cover 
        the previous cone; 
        successful regeneration then occurs in the third attempt at time $t_3$.}
      \label{fig:MBNG-cones}
\end{figure}

With these ingredients, we can discuss the regeneration construction: 
Assume that we find a space-time `cone' $C$ (with fixed suitable base diameter and slope) 
which is centred at the current space-time position $(X_n,-n)$ of the walk such that 
\begin{itemize} 
\item[(i)] 
  $C$ covers the path and everything it has `explored' until the $n$-th step  
  (since the last regeneration),

\item[(ii)] the configuration in $\eta^{\rm stat}$ at the base of the 
  cone $C$ is `good' and

\item[(iii)] `strong' coupling for $\eta^{\rm stat}$ as defined in property~(c) above 
  occurs inside the cone $C$.
\end{itemize}
Then, the conditional law of the future path increments is completely determined 
by the configuration $\eta^{\rm stat}$ at the base of the cone and we can `start afresh.' 
It may happen that in order to find a cone with properties (i)--(iii), several attempts 
are needed, see Figure~\ref{fig:MBNG-cones} for an illustration.
\index{regeneration construction}%

This construction expresses the path increments between the regeneration times as a functional of a 
well-behaved Markov chain (which keeps track of the local configuration at the 
base of the corresponding cones at the regeneration times). Given this, \eqref{MBNG-eq:lbrwLLN} 
and \eqref{MBNG-eq:lbrwaCLT} are fairly standard.

    






  
\smallskip

In ongoing work with Andrej Depperschmidt and Timo Schl\"{u}ter 
we consider the joint dynamics of several ancestral lineages in the logistic branching random walk and 
establish properties analogous to those in Section~\ref{MBNG-sect:DRWOPC} 
for walks on the oriented percolation cluster.

\section{Discussion}
\label{MBNG-sect:notes}

Our ancestral walks with dynamics as in \eqref{MBNG-eq:defXdynamics}, \eqref{MBNG-eq:11}, \eqref{MBNG-eq:lbrwXtransprob} 
are generally speaking random walks in dynamical random environments (RWDRE). 
This is currently a very active field of research and we do not attempt
to give an overview here, but refer to \cite{MBNG-A10} for a good overview of the area up to 2010. 
There are recent papers on random walks in dynamical random conductances, random walks on dynamical percolation, random walks in dynamical random environments given by interacting particle systems as for instance exclusion processes. The general results have often strong assumptions on the environment (mixing conditions, spectral gap assumptions, uniform lower bounds for the transition probabilities of the walk). On the other hand, the `case studies' often refer to specific models and do not provide a general approach.
Hence, this is an area where there is still a lot to understand. 
See, e.g., the recent works \cite{MBNG-ABF18, {MBNG-BHT18}} and the discussion and references there. 
Let us point out that our walks \eqref{MBNG-eq:defXdynamics}, \eqref{MBNG-eq:11}, \eqref{MBNG-eq:lbrwXtransprob} are
somewhat special inside the general class of RWDRE in that the natural `forwards' time direction for the 
walk is `backwards' in time for the environment, whereas often researchers in RWDRE study walks on certain interacting 
particle systems where the walk and the underlying system have the same forwards time direction. Also, let us mention that while in recent work, see \cite{MBNG-BR18}, the assumption of ellipticity of the environment, i.e. on uniform lower bounds for the transition probabilities of the walk, is not present anymore, our model still does not fit in, since our environment is not stationary.
\index{random walk!in dynamic random environment}%
\medskip

\noindent {\bf Acknowledgements.} The authors thank Ji\v{r}\'{i}
\v{C}ern\'{y}, Andrej Depperschmidt, Katja Miller and Sebastian Steiber for the many enlightening
discussions we had in the course of this project. We would also like
to thank Iulia Dahmer, Frederik Klement and Timo Schl\"{u}ter and an
anonymous referee for carefully reading the manuscript and for their
helpful comments.

\printindex


\begin{thebibliography}{99}\itemsep=2pt

\bibitem{MBNG-A10} L.~Avena, \textit{Random walks in dynamic random environments}, 
  Proefschrift Universiteit Leiden (PhD dissertation), 2010.
  \url{http://hdl.handle.net/1887/16072}
  
\bibitem{MBNG-ABF18} L.~Avena, O. Blondel, A. Faggionato, 
Analysis of random walks in dynamic random environments via $L^2$-perturbations,
\textit{Stoch.\ Proc.\ Appl.} \textbf{128} (2018), 
3490--3530.

\bibitem{MBNG-BartonDepaulisEtheridge:2002}
N.~H.~Barton, F.~Depaulis, and A.~M.~Etheridge, 
Neutral evolution in spatially continuous populations.
\textit{Theoret. Popul. Biol.} {\bf 61} (2002), 
31 -- 48.

\bibitem{MBNG-BSW18}
G.~Berzunza, A.~Sturm, A.~Winter,
Trait-dependent branching particle systems with competition and multiple offspring,
preprint, \texttt{arXiv:1808.09345}.

\bibitem{MBNG-BD07}
M.~Birkner and A.~Depperschmidt, 
Survival and complete convergence for a spatial branching system with
local regulation. \textit{Ann. Appl. Probab.} \textbf{17} 
(2007), 1777--1807.

\bibitem{MBNG-BCDG13} 
M.~Birkner, J.~{\v{C}}ern{\'y}, A.~Depperschmidt and N.~Gantert, 
Directed random walk on the backbone of an oriented percolation cluster, 
\textit{Electron. J. Probab.} \textbf{18} (2013), no.~80, 35p.

\bibitem{MBNG-BCD16} 
M.~Birkner, J.~{\v{C}}ern{\'y} and A.~Depperschmidt, 
Random walks in dynamic random environments and ancestry under local population regulation, 
\textit{Electron. J. Probab.} \textbf{21} (2016), no.~38, 1--43.

\bibitem{MBNG-BGS18} 
M.~Birkner, N.~Gantert and S.~Steiber, 
Coalescing directed random walks on the backbone of a 1+1-dimensional 
oriented percolation cluster converge to the Brownian web,
\textit{ALEA Lat. Am. J. Probab. Math. Stat.} \textbf{16} (2019), 1029--1054.

\bibitem{MBNG-BS19} 
M.~Birkner, R.~Sun, 
Low-dimensional lonely branching random walks die out, 
\textit{Ann. Probab. } \textbf{47} (2019), 
774--803.

\bibitem{MBNG-BR18}
M.~Biskup, P.-F.~Rodriguez,
Limit theory for random walks in degenerate time-dependent random environments,
\textit{Journal of Functional Analysis} \textbf{274} (2018), 
985--1046.

\bibitem{MBNG-BlathEtheridgeMeredith:2007}
J.~Blath, A.~M.~Etheridge, and M.~Meredith,
Coexistence in locally regulated competing populations and survival
  of branching annihilating random walk,
  \textit{Ann. Appl. Probab. } \textbf{17} (2007), 
  1474--1507.

\bibitem{MBNG-BlathHammerOrtgiese:2016}
J.~Blath, M.~Hammer, M.~Ortgiese,
The scaling limit of the interface of the continuous-space symbiotic branching model,
\textit{Ann. Probab.} \textbf{44} (2016), 
807--866.


\bibitem{MBNG-BHT18} O.~Blondel, M.~R.~Hil\'{a}rio, A.~Teixeira, 
Random Walks on Dynamical Random Environments with Non-Uniform Mixing, 
preprint, \texttt{arXiv:1805.09750}.

\bibitem{MBNG-BolkerPacala:97}
 B.~M.~Bolker, S.~W.~Pacala 
Using moment equations to understand stochastically driven spatial
  pattern formation in ecological systems,
\textit{Theoret. Popul. Biol.} \textbf{52}, (1997) 179--197.

\bibitem{MBNG-BolkerPacala:99}
 B.~M.~Bolker, S.~W.~Pacala,
Spatial moment equations for plant competition: Understanding spatial
  strategies and the advantages of short dispersal,
\textit{American Naturalist} \textbf{153}, (1999) 575--602.

\bibitem{MBNG-BolthausenSznitman:2002}
E.~Bolthausen and A.-S.~Sznitman, {On the static and dynamic points of
  view for certain random walks in random environment}, \textit{Methods Appl. Anal.
  \textbf{9}} (2002), 
345--375, Special issue dedicated to Daniel W.
  Stroock and Srinivasa S. R. Varadhan on the occasion of their 60th birthday.


\bibitem{MBNG-Durrett:84}
R.~Durrett, {Oriented percolation in two dimensions}, \textit{Ann. Probab.}
\textbf{12} (1984), 
999--1040. 

\bibitem{MBNG-Etheridge:2004}
 A.~M.~Etheridge,
Survival and extinction in a locally regulated population.
\textit{Ann. Appl. Probab.} \textbf{14} 
(2004), 188--214.

\bibitem{MBNG-Etheridge:2011} 
A.~M.~Etheridge,
\textit{Some mathematical models from population genetics}, Volume 2012
  of \textit{Lecture Notes in Mathematics}, Springer, 2011. 
Lectures from the 39th Probability Summer School held in Saint-Flour, 2009.

\bibitem{MBNG-FournierMeleard:2004}
N.~Fournier, S.~M{\'e}l{\'e}ard,
A microscopic probabilistic description of a locally regulated
  population and macroscopic approximations.
  \textit{Ann. Appl. Probab.} \textbf{14} 
  (2004), 1880--1919.

\bibitem{MBNG-GrevenSturmWinterZaehle:19+}
A.~Greven, A.~Sturm, A.~Winter and I.~Z\"{a}hle,
Multi-type spatial branching models for local self-regulation I: 
Construction and an exponential duality, preprint, \texttt{arXiv:1509.04023}.  


\bibitem{MBNG-GrimmetHiemer:02}
G.~Grimmett and P.~Hiemer, \textit{Directed percolation and random walk}, In and
  out of equilibrium ({M}ambucaba, 2000), Progr. Probab., vol.~51, Birkh\"auser, 2002, pp.~273--297.

\bibitem{MBNG-GKS:13}
 O.~G\"{u}n, W.~K\"{o}nig, O.~Sekulovi\'{c}, 
Moment asymptotics for branching random walks in random environment. 
\textit{Electron. J. Probab.} \textbf{18} (2013), no.~63, 1--18.

\bibitem{MBNG-GKS:15}
O.~G\"{u}n, W.~K\"{o}nig, O.~Sekulovi\'{c}, 
Moment asymptotics for multitype branching random walks in random environment.
\textit{J. Theoret. Probab.} \textbf{28} (2015), 
1726--1742.

\bibitem{MBNG-HammerOrtgieseVoellering:18}
M.~Hammer, M.~Ortgiese, F.~V\"{o}llering, A new look at duality for the symbiotic branching model, 
\textit{Ann. Probab.} \textbf{46} (2018), 
2800--2862.

\bibitem{MBNG-HammerOrtgieseVoellering:19+}
M.~Hammer, M.~Ortgiese, F.~V\"{o}llering, 
Entrance laws for annihilating Brownian motions and the continuous-space voter model,
preprint, \texttt{arXiv:1801.06197}.

\bibitem{MBNG-HutzenthalerWakolbinger:07}
M.~Hutzenthaler, A.~Wakolbinger,
Ergodic behaviour of locally regulated branching populations.
\textit{Ann. Appl. Probab.} \textbf{17} 
(2007), 474--501.

\bibitem{MBNG-Kallenberg} 
O.~Kallenberg, Stability of critical cluster fields.
\textit{Math. Nachr.} \textbf{77} (1977), 7--43.

\bibitem{MBNG-Kuczek:1989} T.~Kuczek, The central limit theorem for
  the right edge of supercritical oriented percolation,
  \textit{Ann. Probab.} \textbf{17} (1989), 
  1322--1332.

\bibitem{MBNG-LawDieckmann:02}
R.~Law, U.~Dieckmann 
Moment approximations of individual-based models.
\textit{In U.~Dieckmann, R.~Law, and J.~A. Metz (Eds.)} {The Geometry of
  Ecological Interactions}, Cambridge Univ. Press, 2002, pp.~252--270.

\bibitem{MBNG-LePardouxWakolbinger:2013}
V.~Le, E.~Pardoux, and A.~Wakolbinger,
``{T}rees under attack'': a {R}ay-{K}night representation of
  {F}eller's branching diffusion with logistic growth,
  \textit{Probab. Theory Related Fields} \textbf{155} 
  (2013), 583--619.

\bibitem{MBNG-M16}
K.~Miller,
Random walks on weighted, oriented percolation clusters, 
\textit{ALEA Lat. Am. J. Probab. Math. Stat.} \textbf{13} (2016), 53--77. 
Erratum: \textit{ALEA Lat. Am. J. Probab. Math. Stat.} \textbf{14} (2017), 173--175. 

\bibitem{MBNG-M17}
  K.~Miller, Random walks on oriented percolation and in recurrent environments,
  Dissertation, Technische Universit\"at M\"unchen, 2017.
  \url{http://mediatum.ub.tum.de/670560?show\_id=1366085}

\bibitem{MBNG-Neuhauser:1992}
C.~Neuhauser, Ergodic theorems for the multitype contact process,
\textit{Probab. Theory Related Fields} \textbf{91} (1992), 467--506.

\bibitem{MBNG-NeuhauserPacala:1999}
C.~Neuhauser, S.~W.~Pacala 
An explicitly spatial version of the {L}otka-{V}olterra model with
  interspecific competition.
  \textit{Ann. Appl. Probab.} \textbf{9} 
  (1999), 1226--1259. 

\bibitem{MBNG-NewmanRavishankarSun2005}
C.~M.~Newman, K.~Ravishankar and R.~Sun,
  {Convergence of coalescing nonsimple random walks to the {B}rownian
    web}, \textit{Electron. J. Probab.} \textbf{10} (2005), no. 2, 21--60.

\bibitem{MBNG-RaghibHillDieckmann:2011}
M.~Raghib, N.~A.~Hill, and U.~Dieckmann, 
A multiscale maximum entropy moment closure for locally regulated
  space time point process models of population dynamics,
  \textit{J. Math. Biol.} \textbf{62} 
  (2011), 605--653.


\bibitem{MBNG-Sawyer:1976} 
S.~Sawyer, 
Results for the stepping stone model for migration in population genetics,
\textit{Ann. Probability} \textbf{4} 
(1976), 699--728.

\bibitem{MBNG-SSS17}
  E.~Schertzer, R.~Sun and J.M.~Swart, 
  {The {B}rownian web, the {B}rownian net, and their universality},
  \textit{Advances in disordered systems, random processes and some applications},
  Cambridge Univ. Press, 2017, pp.~270--368.

\bibitem{MBNG-S17}
  S.~Steiber, Ancestral lineages in the contact process: scaling and hitting properties, 
  Dissertation, Johannes Gutenberg-Universit\"{a}t Mainz, 2017.
  \url{https://nbn-resolving.org/urn:nbn:de:hebis:77-diss-1000010573}

\bibitem{MBNG-VeberWakolbinger:2014}
A.~V{\'e}ber, A.~Wakolbinger,
The spatial {L}ambda-{F}leming--{V}iot process: {A}n event-based
  construction and a lookdown representation,
  \textit{Ann. Inst. Henri Poincar\'e Probab. Stat.} \textbf{51} 
  (2015), 570--598.

\bibitem{MBNG-WeissKimura:1965}
G.~H.~Weiss, and M.~Kimura, A mathematical analysis of the stepping stone model of genetic
  correlation,
\textit{J. Appl. Probability} \textbf{2} (1965), 129--149.

\bibitem{MBNG-Wilkinson-Herbots:2003}
H.~M.~Wilkinson-Herbots, Coalescence times and {$F_{ST}$} values in subdivided populations
with symmetric structure,
\textit{Adv. in Appl. Probab.} \textbf{35} 
(2003), 665--690.

\bibitem{MBNG-ZaehleCoxDurrett:2005}
I.~Z{\"a}hle, , J.~T.~Cox and R.~Durrett (2005).
The stepping stone model, II. {G}enealogies and the infinite sites model,
\textit{Ann. Appl. Probab.} {\bf 15} 
(2005), 671--699.

\end{thebibliography}
\end{document}